\documentclass[letterpaper,10pt]{article}

\usepackage[draft]{hyperref}
\usepackage{amsmath}
\usepackage{amssymb}
\usepackage{psfrag}
\usepackage{graphicx}
\usepackage{here}
\let\epsilon =\varepsilon

\newtheorem{Thm}{Theorem}
\newcommand{\tn}{\textnormal}

\title{Robust control of a bimorph mirror\\ for adaptive optics system}

\date{}

\author{Lucie Baudouin, Christophe Prieur, \\
Fabien Guignard and Denis Arzelier\footnote{e-mail: {\tt baudouin@laas.fr, arzelier@laas.fr, cprieur@laas.fr, fguignard@laas.fr}}\\
{\it\footnotesize LAAS - CNRS; Universit\'e de Toulouse; 7, avenue du Colonel Roche, F-31077 Toulouse, France.}}

\begin{document}

\maketitle

\begin{abstract}
We apply robust control technics to an adaptive optics system including a dynamic
model of the deformable mirror. The dynamic model of the mirror is a modification of the usual plate
 equation. We propose also a state-space approach to model the turbulent phase.
 A continuous time control of our model is suggested taking into account the frequential behavior
 of the turbulent phase. An $H_\infty$ controller is designed in an infinite dimensional setting.
 Due to the multivariable nature of the control problem involved in adaptive optics systems,
 a significant improvement is obtained with respect to traditional single input single output methods.\\
\end{abstract}


\noindent \textit{OCIS:}  010.1330, 120.4640.\\

%%%%%%%%%%%%%%%%%%%%%%%%%%
\section{Introduction}\label{1}
%%%%%%%%%%%%%%%%%%%%%%%%%%

For several decades it has been now possible to use adaptive optic (AO) systems to actively correct the  distortions affecting an incident wavefront propagating through a turbulent medium. A particularly interesting application of this technique is in the field of astronomical ground-based
imaging. The idea behind AO systems is to generate a corrected
wavefront as close as possible to the genuine incident plane
wavefront thanks to a deformable mirror (DM). An AO system
is also composed of a wavefront sensor measuring the resulting
distortion of the collected wavefront after correction by the DM. Based on these measured
signals, the voltage applied to the piezoelectric actuators is computed in order to reshape the mirror.
The tilts (first order modes) of the wavefront are corrected by a first mirror.
Then, the DM is part of the control-loop for the correction of higher-order modes of the wave front.
Different types of sensors (curvature sensor, pyramid wavefront sensor) may be used to estimate the distortions affecting the incoming wave-front but the most common encountered
in existing applications is the Shack-Hartmann (SH) sensor.
There also exists different type of deformable mirrors
and we choose to study the case of the most common one.
For additional details on basic principles of adaptive optics, see \cite{R}, e.g..

\begin{figure}[t]
\centerline{\includegraphics[width=10cm]{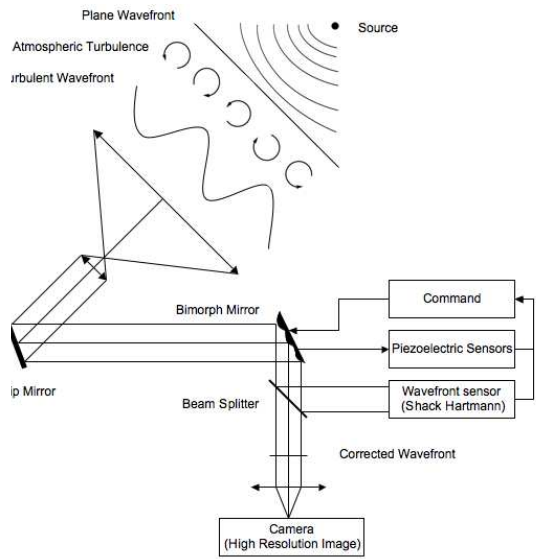}}
\caption{An adaptive optics system -
\footnotesize{The control loop consists in a SH sensor analyzing the incoming
wavefront
and a layer of piezoelectric sensors giving the precise position of the bimorph mirror,
both of them allowing the calculation of the appropriate command of the DM
in order to recover the genuine wavefront.}}\label{Principe}
\end{figure}

This paper is devoted to the design of control laws for
an adaptive optics system formed by a bimorph mirror and a
Shack-Hartmann sensor (see Figure \ref{Principe}). Most often, the
existing adaptive optics systems use static models and very
basic control algorithms based on frequent measurements of the influence of
each actuator of the mirror to each output of the SH. This allows the computation of
an interaction matrix gathering the corresponding influence
functions. Here, our goal is to consider the design of an adaptive optics system control loop
from a modern automatic control point of view as in \cite{RAY06}  and \cite{F-T}. This
means first that dynamics of the different elements involved in
the control-loop have to be taken into account. In particular, a
specific dynamic model for the DM is proposed for control purpose (as already presented in \cite{B-P-A} and see also \cite{M-G}). Secondly, a state-space model of the turbulent phase,
built from its frequency domain characteristics, is defined \cite{CRM95}.

The main contribution concerns the infinite dimension setting introduced in this
paper.
More precisely, while in the literature, only static finite dimensional models are
considered,
a model based on a particular partial differential equation (pde) is used for the DM.
We believe that our point of view matches well with the reduction of the size of the
actuators
and the significant augmentation of their numbers in many devices, as in AO for Very
Large Telescopes.

In reference \cite{L-P}, a thin elastic plate model of a deformable bimorph mirror
is derived.
This model is based on a periodic distribution of embedded
piezoelectric patches that may be used as sensors or actuators.
The idea is then to elaborate a robust control strategy based on
modern control tools for distributed parameter systems \cite{VK1}. Moreover, in contrast
to \cite{PAS93} and \cite{RAY06}, we do not need to compute any interaction matrix modelling the relation between the input on the piezoelectric patches
attached to the mirror and the output given by the Shack-Hartmann sensor. The interaction matrix can be
seen as a static model of the mirror whereas a more general dynamical model of the mirror is used
here.\\

For the sake of clarity of this study, we emphasize here the main informations
about the frame we choose for our modeling of robust control of an AO system.
We consider a continuous time state-space model of an AO loop (as in \cite{M-G}  and
instead of a discrete one in \cite{RAY06} and \cite{PAS93}) and without delay.
In practice an AO system uses discrete wavefront sensing data with inherent temporal
delays and of course it is possible to derive a discrete time extension of our model but it
is not our point here, even if we recognize that the performance will somehow be affected.
Our contribution relies mainly on the new pde model of the DM and we aim at
using the $H_{\infty}$~control theory for infinite dimension setting in order to
recover
at least similar performance as the one of LQG control for a standard model of the
DM (see \cite{PAS93}). One should notice that our model depends only on a few physical parameters
(such as the density, the stiffness... see the \textit{Bimorph mirror model} subsection below for more details), parameters that could be considered as uncertain quantities the control law should take into account. Therefore,we do not need either to compute an interaction matrix (which is more and more complicate to compute when the number of the sensors and of the actuators increases as for Very Large Telescopes), or the inverse of this interaction matrix \cite{R}.\\

The control problem is solved using an $H_{\infty}$~control setting.
The first motivation is that $H_{\infty}$~control theory provides intrinsic properties of robustness while optimizing on the worst-case performance. Another
motivation is the multivariable nature of the control problem
involved in adaptive optics system design \cite{F-T}.
Current adaptive optics control systems use decoupling modal
control to rewrite the original problem as several decoupled single
input single output control problems. Because $H_{\infty}$~control framework
may easily handle a multivariable dynamic model
of the bimorph DM in the synthesis process, the
obtained robust controller outperforms usual static
control approaches of the literature.
In addition, the use of Hinfinity controllers induces, in general, some robustness properties of the closed-loop while $H_2$/$LQG$ controllers (privileged in general, see \cite{PAS93}) lead to improvement of the performance but with no robustness guarantee (see \cite{SP}). So far, we do not claim to have solved the complete problem of AOS synthesis (with delays and limitations of performance introduced by sampling) but we think that this new setting will probably address fundamental issues encountered in the very large telescopes context. This work is meant to illustrate the realizability of such an approach on realistic instances of AOS Design.\\

The outline of the paper is the following. First, the adaptive
optics control system is described (see Section \ref{2}) through the presentation of
the models of the bimorph mirror and the turbulent phase.
The third section is dedicated to the robust $H_{\infty}$~control
setting in the infinite dimension framework and its formulation in our
particular case. The last section contains the description of the truncated
model and the numerical results.

%%%%%%%%%%%%%%%%%%%%%%%%%%%
\section{The adaptive optics model}\label{2}
%%%%%%%%%%%%%%%%%%%%%%%%%%%

The bimorph mirror is
composed of a purely elastic and reflective plate equipped with piezoelectric actuators (in
order to deform the shape of the mirror) and piezoelectric sensors
(to measure the effective deformation). A Shack-Hartmann sensor
then analyzes the resulting phase $\phi_\tn{res}$ of the
wavefront, after reflection in the deformable mirror of the
turbulent phase $\phi_\tn{tur}$.\\

Different types of disturbances have
to be faced with: $w_\tn{mod}$ represents unstructured uncertainty
(neglected dynamics) affecting the model, $w_\tn{piezo}$ and $w_{_\tn{SH}}$
 are noise signals respectively attached to
piezoelectric and Shack-Hartmann sensors. Finally, $\phi_\tn{tur}$ is
the turbulent phase of the wavefront introduced by the atmospheric
perturbation.\\

We denote by $e=e(r,\theta,t)$ the transverse displacement of the circular mirror at point of polar
coordinates $(r,\theta)$ and time $t$, while $\lambda$ is the light wavelength.
The corrected phase produced by $e$ is then
given by $~\phi_\tn{cor}=\frac{4\pi}{\lambda}e~$ leading to a
resulting phase:
\begin{equation}
\label{tobecontrolledopticoutput}
\phi_\tn{res} = -
\frac{4\pi}{\lambda}e + \phi_\tn{tur}
\end{equation}
The optic sensor's output, computed by Shack-Hartmann sensor is:
\begin{equation}
\label{opticoutput}
y_{_\tn{SH}}=- \frac{4\pi}{\lambda} e + \phi_\tn{tur} + c w_{_\tn{SH}}
\end{equation}
where $c$ is a modelling parameter of the perturbation.\\

Finally, we note that the control input is the voltage $u$
applied to the piezoelectric actuators and the corresponding
piezoelectric output is the voltage $y_{pe}$ measured with the
piezoelectric inclusions used as sensors (see equations (\ref{plate}) and (\ref{piezooutput}) below).
Indeed, in comparison with many other devices, where the only information used to compute the voltage
$u$ comes from the wavefront analyzer, the additional possibility of measuring the deflection of the mirror through a layer of piezoelectric sensors (see Figure \ref{Principe}) is considered here.

It is recalled that the goal of the adaptive optics control system is to minimize the resulting phase of the wavefront
using Shack-Hartmann measurements.

%%%%%%%%%%%%%%%%%%%%%%%%%%%%
\subsection*{Bimorph mirror model}\label{22}
%%%%%%%%%%%%%%%%%%%%%%%%%%%%

To obtain the model of a bimorph mirror (see an outline in \cite{B-P-A}), we consider three different layers. One is purely elastic and reflective,
the second one is equipped with piezoelectric inclusions used as
actuators, the third one is equipped with piezoelectric inclusions
used as sensors.  The heterogeneities are periodically distributed.  In reference \cite{L-P},
the authors derive the following dynamical model of the mirror (a partial differential equation with respect to $(r,\theta,t)$):
\begin{equation}
\label{plate} \rho~\partial_{tt} e+Q_{1}\Delta^2 e
+Q_{2}e =\widetilde{d}_{31}\Delta u + \rho b w_\tn{mod}
\end{equation}
with the initial conditions
$e(r,\theta,t=0)=e_0(r,\theta)$ and $\partial_t e (r,\theta,t=0)= e_1(r,\theta)$.
The voltage $y_{pe}$ computed by the piezoelectric sensors
is given by
\begin{equation}\label{piezooutput}
y_{pe} = \tilde{e} _{31} \Delta e+dw_\tn{pe}.
\end{equation}
The following notations are defined:
\begin{itemize}
  \item $(r,\theta)$ are the spatial coordinates of a point of the disk $\Omega$ of
radius $a$ and $t$ is the time;
  \item $\Delta$ is the Laplacian operator and for a general function $v(r,\theta)$ in polar coordinates
   $$\Delta v = \dfrac{\partial v}{\partial
r^2}+\dfrac{1}{r}\dfrac{\partial v}{\partial r}+\dfrac{1}{r^2}\dfrac{\partial v}{\partial \theta ^2};$$
  \item $u$ is the voltage applied to the inclusions of the actuator layer;
  \item $\rho$ is the surface density, $\nu$ is the Poisson ratio of the mirror's material, $Q_1$ is the stiffness coefficient
  and $Q_2$ is a correction coefficient;
  \item   $\tilde{e}_{31}$ and $\widetilde d_{31}$
  are proportional to the piezoelectric tensor coefficient $d_{31}$ (for more physical details see \cite{N});
  \item $b$ and $d$ are linear applications on appropriate spaces;
  \item $w_\tn{mod}$ and $w_\tn{pe}$ are unknown perturbations
modelling the model errors of the plate equation and the measurement noise of the piezoelectric output.
\end{itemize}

The boundary conditions are those of the free edges case (VLT and the experimental device SESAME,
see Subsection~\ref{numerical}):
\begin{equation}\label{bound}
\begin{array}{ll}
\left.\dfrac{\partial ^2 e}{\partial r^2}+\nu\left(\dfrac{1}{r}\dfrac{\partial e}{\partial r}
+\dfrac{1}{r^2}\dfrac{\partial^2 e}{\partial \theta ^2}\right)\right|_{r=a}&=0\vspace{0,2cm}\\
\left.\dfrac{\partial}{\partial r}\left(\Delta e\right)+\dfrac{1}{r}(1-\nu)\dfrac{\partial}{\partial r}
\left(\dfrac{1}{r}\dfrac{\partial e}{\partial \theta}\right)\right|_{r=a}&=0
\end{array}
\end{equation}

%%%%%%%%%%%%%%%%%%%%%%
\subsection*{Turbulent phase model}\label{23}
%%%%%%%%%%%%%%%%%%%%%%

In order to complete our optics system model, we need to develop a model of the turbulence phase.

A usual representation of atmospheric phase distortion is made through the orthogonal basis of Zernike polynomials because the first Zernike modes
correspond to the main optical aberrations. An infinite number of Zernike functions is required to characterize the wavefront, but a truncated basis
is used in general for implementation purpose. Note that a 14-th order approximation contains 92\% of the phase information, without taking into
account the piston mode which represents the average phase distortion \cite{PAS93}. The tip/tilt modes are not part of our modelling of the turbulent
phase because of their correction by a dedicated mirror.
We will therefore work with the $12$ first modes of Zernike given in reference \cite{Nol76} and recalled here (see Table \ref{Zernike}), excluding the three first ones.

\begin{table}[h]
\begin{center}
\begin{tabular}{|c c c|c|c|}
\hline
$i$&\(n\)&\(m\)&\(Z_i(r,\theta)\)%&Name
\\
\hline
1&0&0&1%&Piston
\\
2&1&1&2\(\frac{r}{a}\cos{\theta}\)%&y tilt
\\
3&1&1&2\(\frac{r}{a}\sin{\theta}\)%&x tilt
\\
4&2&0&\(\sqrt{3}(2(\frac{r}{a})^2-1)\)%&Focus
\\
5&2&2&\(\sqrt{6}(\frac{r}{a})^2\cos{2\theta}\)%&Astigmatism
\\
6&2&2&\(\sqrt{6}(\frac{r}{a})^2\sin{2\theta}\)%&Astigmatism
\\
7&3&1&\(\sqrt{8}(3(\frac{r}{a})^3-2\frac{r}{a})\cos{\theta}\)%&Pure coma
\\
8&3&1&\(\sqrt{8}(3(\frac{r}{a})^3-2\frac{r}{a})\sin{\theta}\)%&Pure coma
\\
9&4&0&\(\sqrt{5}(6(\frac{r}{a})^4-6(\frac{r}{a})^2+1)\)%&Spherical
\\
10&3&3&\(\sqrt{8}(\frac{r}{a})^3\cos{3\theta}\)%&Trefoil
\\
11&3&3&\(\sqrt{8}(\frac{r}{a})^3\sin{3\theta}\)%&Trefoil
\\
12&4&2&\(\sqrt{10}(4(\frac{r}{a})^4-3(\frac{r}{a})^2)\cos{2\theta}\)%&Fifth-order astigmatism
\\
13&4&2&\(\sqrt{10}(4(\frac{r}{a})^4-3(\frac{r}{a})^2)\sin{2\theta}\)%&Fifth-order astigmatism
\\
14&4&4&\(\sqrt{12}(10(\frac{r}{a})^5-12(\frac{r}{a})^3+3(\frac{r}{a}))\cos{4\theta}\)
\\
15&4&4&\(\sqrt{12}(10(\frac{r}{a})^5-12(\frac{r}{a})^3+3(\frac{r}{a}))\sin{4\theta}\)
\\
\hline
\end{tabular}
\caption{ First 15 Zernike Functions}
\label{Zernike}
\end{center}
\end{table}

The turbulent phase $\phi_\tn{tur}$ is approximated as follows:
\[ \phi_\tn{tur}(r,\theta,t) \approx \sum_{i=4}^{N_Z}\phi_i(t)Z_i(r,\theta) \]
where \(N_Z \geq 15\). \(Z_i\) is the \(i\)-th Zernike function and for all $i$, \(\phi_i(t)\) is a random time-varying coefficient corresponding to
the projection of \(\phi_\tn{tur}\) on \(Z_i\).

\begin{figure}[h]
\psfrag{a}[][][0.6]{$w=\left(\begin{array}{c} w_1(j\omega) \\ \vdots \\ w_{N_Z}(j\omega) \end{array}\right)$}
\psfrag{c}[][][0.6]{$\left(\begin{array}{c} \phi_1(j\omega)\\ \vdots \\ \phi_{N_Z}(j\omega) \end{array}\right)=\phi$}
\psfrag{e}[][][1]{$\huge{H(jw)}$} \centerline{\includegraphics[width=6cm]{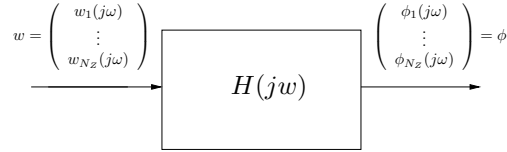}}
\caption{Shaping filter generating $\phi$ -
\footnotesize{The turbulent phase $\phi$ is modeled through a linear shaping filter of transfer function $H$ from the noise $w$  }}
\label{schema_principe}
\end{figure}

To build a state-space representation of the turbulent phase, \(\phi_\tn{tur}\) is modelled as the output of a linear shaping filter (illustrated by Figure \ref{schema_principe}) of the form :
\begin{equation}
\label{eqtur}
\phi '=F\phi+Gw
\end{equation}
where $\phi=(\phi_4,\cdots,\phi_{N_Z})$, $w=(w_4,\cdots,w_{N_Z})$,  \(F\) and \(G\) are two time-invariant square matrices of $(N_Z-3)$-dimension and
\(w\) is a stationary zero-mean white gaussian noise. \(\phi_\tn{tur}\) is therefore a stationary process.

In order to compute $F$ and $G$, the results presented in \cite{CRM95} and based on the Kolmogorov theory of turbulence and associated approximations
in the frequency domain are used here. They confirm similar results proposed in \cite{HB76} and complete the study of frequency domain behavior for
each Zernike coefficient. Each Zernike function's spectrum are characterized by a cut-off frequency whose heuristic expression is given by:
\begin{equation}
f_{c_i} \sim 0.3(n_i+1)\dfrac{V}{D}
\end{equation}
where $n_i$ is the radial order of the Zernike number $i$, $V$ is the average wind-speed and $D$ the diameter of the circular aperture of the telescope.

The random process $\phi$ is supposed to be composed of $N_Z-3$ decoupled first-order Markov processes. For $i=4 \cdots N_Z$, we have:
\begin{equation}
H_i(P)=\dfrac{\phi_i(p)}{\omega_i(p)}=\dfrac{1}{1+\tau_ip} ~\textnormal{with }~ \tau_i=\dfrac{1}{2\pi f_{c_i}}
\end{equation}
In other words, $F=\tn{diag}_i(-\dfrac{1}{\tau_i})$.

\begin{table}[h]
\begin{center}
\begin{tabular}{|c c|c|c||c c|c|c|}
\hline
$i$&$j$&$F_{i,j}$&$G_{i,j}$&$i$&$j$&$F_{i,j}$&$G_{i,j}$\\
\hline
1&1&-508,9&27.10	&6&6&-848.2&10.20	\\
1&6&0&-4.499		&7&7&-678.6&16.38	\\
2&2&-508,9&27.11	&8&8&-678.6&16.38	\\
2&9&0&-4.455		&9&2&0&-4.455 \\
3&3&-508,9&27.11	&9&9&-848.2&10.64 \\
3&10&0&-4.455	&10&3&0&-4.455		\\
4&4&-678.6&15.48	&10&10&-848.2&10.64 \\
4&11&0&-3.555	&11&4&0&-3.555 \\
5&5&-678.6&15.47	&11&11&-101.8&8.047 \\
5&12&0&-3.555	&12&5&0&-3.555\\
6&1&0&-4.499		&12&12&-101.8& 8.047\\
\hline
\end{tabular}
\vspace{0.2cm}
\caption{Atmospheric phase distortion state-space model
with the average wind-speed $V=9$m s$^{-1}$, \(\dfrac{D}{r_0}=8\) and the wavelength $\lambda=550$nm} \label{matrix}
\end{center}
\end{table}

The matrix $G$ is obtained from the steady-state Lyapunov equation verified by the correlation matrix $P_{\phi}(\infty)$:
\begin{equation}
GG'=-(FP_{\phi}(\infty)+P_{\phi}(\infty)F')
\end{equation}

A closed-form expression for the spatial covariance matrix is given in \cite{Nol76}.
\begin{eqnarray*}
&&P_{\phi}(\infty)  =\tn{cov}(\phi_i,\phi_j)=E(\phi_i\phi_j)\\
&&=7.19\times
10^{-3}\times(-1)^{(n_i+n_j-m_i-m_j)/2}\left(\frac{D}{r_0}\right)^{\frac{5}{3}}\\
&& \times\sqrt{(n_i+1)(n_j+1)}\pi^{\frac{8}{3}}\\
& &\times\dfrac{\Gamma\big(\dfrac{14}{3}\big)\Gamma\big(\dfrac{n_i+n_j-\frac{5}{3}}{2}\big)}
{\Gamma\big(\dfrac{n_i-n_j+\frac{17}{3}}{2}\big)\Gamma\big(\dfrac{n_i-n_j+\frac{17}{3}}{2}\big)
\Gamma\big(\dfrac{n_i+n_j+\frac{23}{3}}{2}\big)}
\end{eqnarray*}
where \(\Gamma\) is the Gamma function and \(r_0\) is the Fried parameter (corresponding to the strength of the turbulence \cite{R}). Table \ref{matrix} shows the non zero entries of the
matrices \(F\) and $G$ for $V=9 \textnormal{m s}^{-1}$ and \(\dfrac{D}{r_0}=8\) (as in reference \cite{RAY06}).

%%%`%%%%%%%%%%%%%%%
\section{Robust Control Results}\label{3}
%%%%%%%%%%%%%%%%%%

The point of this section is to prove that the new model we propose for AO systems is valid for an $H\sb \infty$-control study. One of the difficulties comes from the infinite dimensional setting. For a survey of the $H\sb \infty$-control theory for the infinite-dimensional case, the interested
reader may have a look at \cite{B-B} or \cite{C-P-VK} for the state-feedback case and \cite{VK1} for the output-feedback case. The main results are a generalization of finite-dimensional regular $H\sb \infty$-control problems (see for instance \cite{SP}). In particular, the solution will be given in terms of the solvability of two coupled Riccati equations.\\

The linear infinite-dimensional model derived from the partial
differential equations presented in Section~\ref{22} has to fit in the following
standard formalism of measurement-feedback control
\begin{equation}\tag{P}\label{P}
\left\lbrace
\begin{array}{ll}
x' =  Ax + B_1 w + B_2 u \\
z  =  C_1 x + D_{12} u \\
y  =  C_2 x + D_{21} w
\end{array} \right.
\end{equation}
where $x$ is the state of the system, $u$ is the control input, $w$ is the
disturbance input, $y$ is the measured output and $z$ is the controlled output.

\begin{figure}[h]
\begin{center}
\includegraphics[width=5cm]{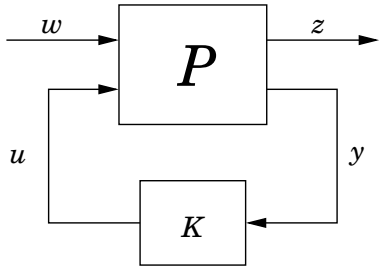}
\end{center}
\caption{Closed-loop system}\label{CLSys}
\end{figure}

Therefore, we introduce the following notations:
\begin{itemize}
\item the state vector $x = (e,\partial_t e,\phi_\tn{tur})^T$ where $e$ is the transverse displacement of the plate and $\phi_\tn{tur}$
is the projection of the turbulent phase on the $N_z$ first Zernike modes;
\item the exogenous disturbance inputs vector $w = (w_\tn{mod}, w_{_\tn{SH}}, w_\tn{tur}, w_\tn{pe})^T$ gathers the different perturbation signals
(uncertainty affecting dynamics of the model and of the turbulence phase, noise vectors of the wavefront analyzer and of piezoelectric sensors);
\item the control inputs vector $u$ is  the voltage applied to piezoelectric patches;
\item the measurement outputs vector
$y =(y_\tn{pe},y_{_\tn{SH}})$ is composed with the piezoelectric
and the wavefront analyzer measured outputs;
\item the controlled outputs vector $z = ( \phi_\tn{res},u)$ contains an optical
part (the resulting phase, see (\ref{tobecontrolledopticoutput})) and
the control input vector~$u$.
\end{itemize}

The aim is to find a dynamic measurement-feedback controller $K$ ensuring that the influence of $w$ on $z$ is smaller than some specific bound. The corresponding standard block diagram is given by Figure \ref{CLSys}.

The controller is assume to have the following form:
\begin{equation}\tag{K}\label{K}
\left\lbrace
\begin{array}{llll}
p' &=&  Mp + Ny\\
u  &=& Lp + Ry
\end{array} \right.
\end{equation}
where $M$ is the infinitesimal generator of a $C_0$-semigroup on a real separable Hilbert space and $N$, $L$ and $R$ are bounded
linear operators. With this controller, the closed-loop system can easily be derived and defines a bounded linear map $S_K$
such that $z(t) = (S_Kw)(t)$. Its bound is denoted $||S_{K}||_{\infty}$.

\begin{figure}[h]
\psfrag{a}[][][0.5]{$\phi_\tn{tur}$}
\psfrag{b}[][][0.5]{$w_\tn{tur}$}
\psfrag{c}[][][0.5]{$w_\tn{mod}$}
\psfrag{d}[][][0.5]{$z_2$}
\psfrag{e}[][][0.5]{$P$}
\psfrag{f}[][][0.5]{$u$}
\psfrag{g}[][][0.5]{$y_{pe}$}
\psfrag{h}[][][0.5]{$y_{SH}$}
\psfrag{i}[][][0.5]{$d$}
\psfrag{j}[][][0.5]{$w_\tn{pe}$}
\psfrag{k}[][][0.5]{$w_{SH}$}
\psfrag{l}[][][0.5]{$I$}
\psfrag{m}[][][0.5]{$z_1$}
\psfrag{n}[][][0.5]{$\phi_\tn{res}$}
\psfrag{o}[][][0.5]{$\phi_\tn{cor}$}
\psfrag{p}[][][0.35]{$\left(\begin{array}{c} e \\ e' \end{array} \right)$}
\psfrag{q}[][][0.5]{$\left [\tilde{e}_{31}\Delta \hspace{3mm} 0 \right ]$}
\psfrag{r}[][][0.6]{$\left [\frac{4\pi}{\lambda} \hspace{3mm} 0 \right ]$}
\psfrag{s}[][][0.5]{$\left [\begin{array}{c} 0 \\ \tilde{d}_{31}\Delta \end{array} \right ]$}
\psfrag{t}[][][1]{$K$}
\psfrag{u}[][][0.5]{$G$}
\psfrag{v}[][][0.4]{$\left [\begin{array}{c} 0\\b \end{array} \right ]$}
\centerline{\includegraphics[width=13cm,height=8cm]{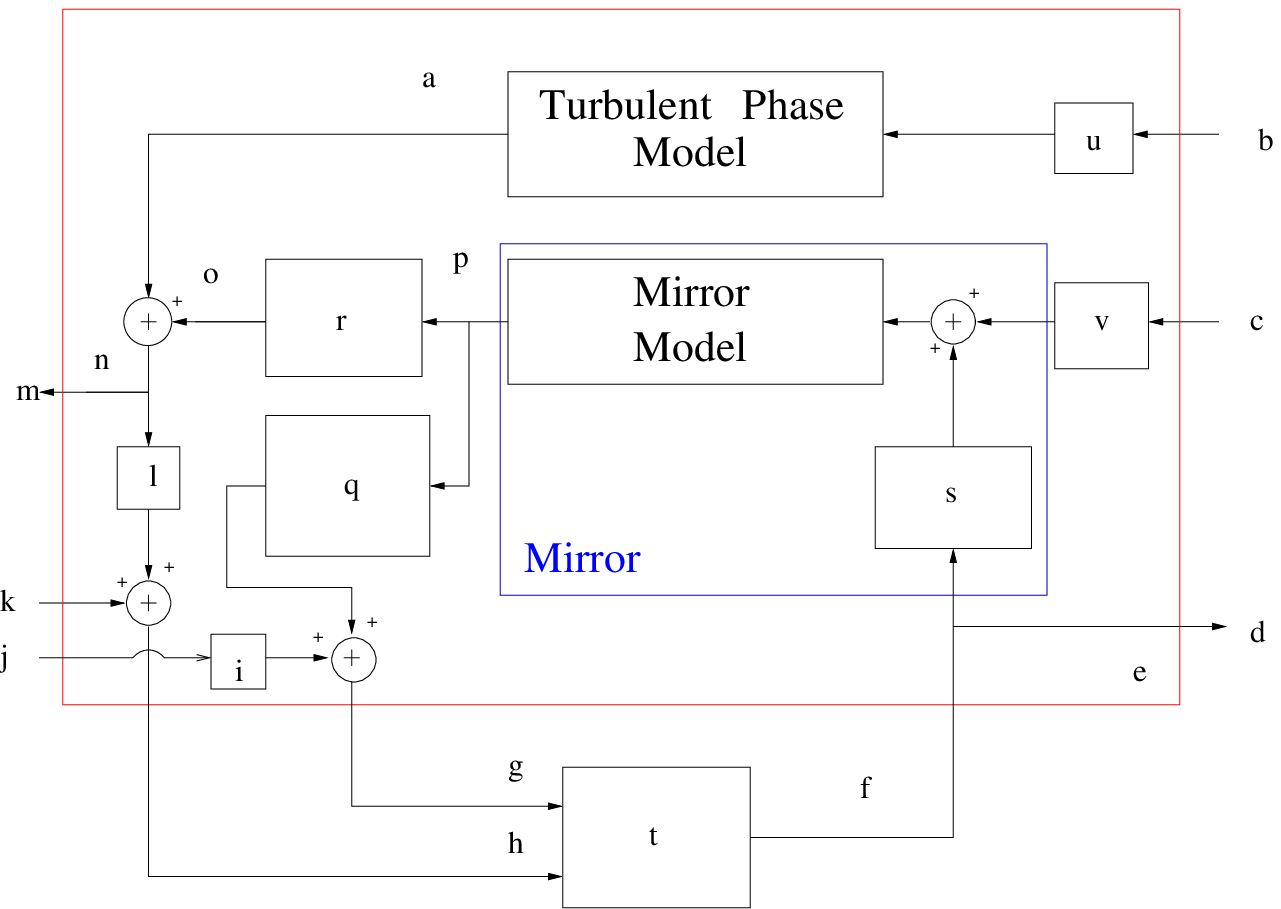}}
\caption{Standard model for adaptive optics system control loop \label{Asserv}}
\end{figure}

The control loop defining the adaptive optics system is sketched in Figure \ref{Asserv}.
If we gather the different equations describing the system, namely (\ref{tobecontrolledopticoutput}), (\ref{opticoutput}), (\ref{plate}), (\ref{piezooutput}) and the forthcoming equation (\ref{phiturinf}) (corresponding to  (\ref{eqtur})), we get
\begin{equation}\label{eqngather}
\left\lbrace\begin{array}{l}
\partial_{tt} e+Q_1\Delta^2 e
+Q_2 e =\widetilde{d}_{31}\Delta u+ b \rho w_\tn{mod} \\
 \partial _{t}\phi_\tn{tur}
= \mathcal F\phi_\tn{tur} + \mathcal G w_\tn{tur} \\
\phi_\tn{res} =\phi_\tn{tur} - \dfrac{4\pi}\lambda e\\
y_{pe} = \tilde{e} _{31} \Delta e+dw_\tn{pe}\\
y_{_\tn{SH}} =\phi_\tn{tur} - \dfrac{4\pi}\lambda e + c w_{_\tn{SH}}.
\end{array}\right.
\end{equation}
Actually, in order to have an unified infinite dimensional modelling of the adaptive optic system's state,
we described the model of $\phi_\tn{tur}$ from equation (\ref{eqtur}) as follows:
\begin{itemize}
\item $\phi_\tn{tur}$ and $w_\tn{tur}$ are the reconstruction of $\phi$ and $w$
on  $N_Z-3$ of the first Zernike modes, such that
$$\phi_\tn{tur} = \sum_{i=4}^{N_Z} \phi_i Z_i~~\textnormal{and}~~w_\tn{tur} = \sum_{i=4}^{N_Z} w_i Z_i$$
\item $\mathcal F $ and $\mathcal G \in \mathcal L(L^2(\Omega))$ satisfy for all $\varphi\in L^2(\Omega)$
\begin{eqnarray*}
\mathcal F (\varphi) &=& \sum_{i=4}^{N_Z} F_{ii}\left<\varphi,Z_i\right>_{L^2(\Omega)} Z_i\\
\mathcal G (\varphi) &=& \sum_{i=4}^{N_Z} \sum_{j=4}^{N_Z}G_{ij}\left< \varphi,Z_j\right>_{L^2(\Omega)} Z_i
\end{eqnarray*}
\end{itemize}
what leads to the $L^2$ turbulent phase model given in (\ref{eqngather})
\begin{equation}\label{phiturinf}
\partial _{t}\phi_\tn{tur}  = \mathcal F\phi_\tn{tur}  + \mathcal G w_\tn{tur}
\end{equation}
where $L^2$ is the Hilbert space of square integrable functions and $\mathcal L(X)$
stands for the set of linear applications on $X$.

Thus, the operators defining the standard form~\ref{P} are built from 
(\ref{eqngather})
$$A = \begin{pmatrix}
  0                  &  I  & 0\\
  -\dfrac{Q_1}{\rho}\Delta^2 -\dfrac{Q_2}{\rho}I  & 0 &  0\\
  0 & 0 & \mathcal F
\end{pmatrix},~~
B_1 = \begin{pmatrix}
0 & 0 & 0 & 0\\
b & 0 & 0 & 0 \\
0 & 0 & \mathcal G & 0
\end{pmatrix},$$

$$B_2 = \begin{pmatrix}
  0     \\
\dfrac{\tilde{d}_{31}  }{\rho} \Delta\\
0
\end{pmatrix},~~
C_1 = \begin{pmatrix}
- \dfrac{4\pi }{\lambda}I  & 0 & I   \\
0                      &  0 & 0
\end{pmatrix},~~
D_{12}  = \begin{pmatrix}
0   \\
I
\end{pmatrix},$$

$$C_2 = \begin{pmatrix}
\tilde{e}_{31}\Delta       & 0        & 0   \\
- \dfrac{4\pi }{\lambda} I &  0 & I
\end{pmatrix},~~
D_{21}  = \begin{pmatrix}
0 & 0 & 0 & d\\
0 & c & 0 & 0
\end{pmatrix}.$$

The appropriate functional spaces associated to the
infinite-dimensional mo\-del are now precisely defined. With the
boundary condition (\ref{bound}), we consider the state space (the mirror $\Omega$ is a disk of radius $a$)
\begin{eqnarray*}
X &=& H^2_\tn{bc}(\Omega)\times L^2(\Omega)\times L^2(\Omega)\\
 &= &\left\{e\in H^2(\Omega), e \textnormal{ satisfying } (\ref{bound}) \right\}
\times \left(L^2(\Omega)\right)^2
\end{eqnarray*}
 the input spaces  $U = H^2(\Omega)\cap H^1_0(\Omega)$ and $W =\left(L^2(\Omega)\right)^4$
 and the output spaces  $Y = Z = \left(L^2(\Omega)\right)^2$,
 where  $H^1_0$ and $H^2$ are the Sobolev spaces
 \begin{gather*}
H^1_0 (\Omega) = \{ f\in L^2(\Omega) /~ \forall i =1,2,~\partial_i f \in L^2(\Omega),~f|_{\partial\Omega}=0  \}\\
H^2 (\Omega) = \{ f\in L^2(\Omega) /~  \forall i,j =1,2,~\partial_i f,~\partial_i\partial_j f \in L^2(\Omega)   \}
\end{gather*}

 This model satisfies all the assumptions of the main theorem of reference  \cite{VK1}.
 We give here a simplified version of this result:

\begin{Thm}\label{VK1}\textnormal{\cite{VK1}}
Let $\gamma>0$. There exists an exponentially
stabilizing dynamic output-feedback controller of the
form $(\ref{K})$ with $\|S_K\|_{\infty} <\gamma$ if and
only if there exist two nonnegative definite operators $P$, $Q\in
\mathcal L(X)$ satisfying the three conditions
\begin{quote}
$(i)$ $\forall x\in D(A)$, $Px\in D(A^*)$,
$$\big(A^*P + PA + P(\gamma^{-2} B_1B_1^* - B_2B_2^*)P + C_1^*C_1\big) x = 0$$
and $A + (\gamma^{-2} B_1B_1^* - B_2B_2^*)P~$ generates an exponentially stable semigroup,\\
$(ii)$ $\forall x\in D(A^*)$, $Px\in D(A)$,
$$\big(AQ + QA^* + Q(\gamma^{-2} C_1^*C_1 - C_2^*C_2)Q + B_1B_1^*\big) x = 0$$
and $A^* + (\gamma^{-2} C_1^*C_1 - C_2^*C_2)Q~$ generates an exponentially stable semigroup,\\
$(iii)$ $$r_\sigma(PQ)<\gamma^2,$$
\end{quote}
where $r_\sigma(PQ)$ stands for the spectral radius of $PQ$.\\
In this case, the controller $K$ given by $(\ref{K})$  and
\begin{equation}
\label{ }
\begin{array}{lll}
M&=& A + (\gamma^{-2} B_1B_1^* - B_2B_2^*)P \\
&&- Q(I-\gamma^{-2} PQ)^{-1}C_2^*C_2\\
N&=& - Q(I-\gamma^{-2} PQ)^{-1}C_2^*\\
L&=& B_2^*P \\
R&=& 0
\end{array}
\end{equation}
is exponentially stabilizing and guarantees that we have
$\|S_K\|_{\infty}<\gamma$, \textit{ie}
$$\|\phi_\tn{res}\|_{L^2(\Omega)} + \|u\|_{L^2(\Omega)} \leq \gamma\|w\|_{\left(L^2(\Omega)\right)^4}.$$
Finally, if the solutions to the Riccati equations exists, then they are unique.
\end{Thm}

Upon additional assumptions that are not detailed here, the main point is to prove that $A$ is the infinitesimal generator of
a $C_0$-semigroup on the real separable Hilbert space $X$. Actually, if we consider the unbounded linear operator
\begin{eqnarray*}
A_1 &:& \mathcal D(A_1) \to X\\
&&\begin{pmatrix} e_0\\ e_1\\ e_2 \end{pmatrix}
\mapsto\begin{pmatrix}
  0                  &  I  & 0\\
  -\Delta^2  & 0 &  0\\
  0 & 0 & 0
\end{pmatrix}
\begin{pmatrix} e_0\\ e_1\\ e_2
\end{pmatrix}
= \begin{pmatrix} e_1\\ -\Delta^2 e_0\\ 0
\end{pmatrix}
\end{eqnarray*}
where
$$ \mathcal D(A_1) =  \left\{e_0\in H^4(\Omega),
e_0 \textnormal{ st } (\ref{bound}) \right\}\times H^2(\Omega)\times L^2(\Omega),$$
then one can prove that $A_1$ is dissipative on $X$. Indeed, we prove that for all $x\in X$,
$$\left< A_1 x , x \right>_{X} \leq 0$$
using  the following scalar product on $H^2_\tn{bc}(\Omega)$ in cartesian coordinates $(x_1,x_2) \in \Omega$,
as suggested in \cite{L-D}:
\begin{gather*}
<u,v>_{H^2_\tn{bc}(\Omega)}\\
=\int_\Omega \Delta u \Delta v - (1-\nu)\left(\dfrac{\partial^2 u}{\partial x_1^2} \dfrac{\partial^2 v}{\partial x_2^2}
+ \dfrac{\partial^2 u}{\partial x_2^2} \dfrac{\partial^2 v}{\partial x_1^2 }\right)\\
+2(1-\nu)\left(\dfrac{\partial^2 u}{\partial x_1\partial x_2} \dfrac{\partial^2 v}{\partial x_1\partial x_2}\right)\,d\Omega. \end{gather*}
Moreover, one can easily check that $A_1$ is also self-adjoint and onto.
Therefore, from Lumer-Phillips' Theorem (see \cite{P}, p. 15), $A_1$
generates a continuous semigroup of linear contractions acting on X.
And finally, since $A$ is the sum of $A_1$ and of a linear
operator bounded on $X$ (as $\mathcal F$ is assumed to be bounded, like $F$),
 the proof is complete (see \cite{G-L-M}, p. 40).\\

 Of course, from a numerical point of view, we need to get an appropriate finite dimensional model.

%%%%%%%%%%%%%%%%%%%
\section{A truncated model for numerical design}\label{33}
%%%%%%%%%%%%%%%%%%%
\subsection{Truncation}
The corresponding finite dimensional model can be presented as :
\begin{equation}\label{Sn}
\left\lbrace
\begin{array}{ll}
x_N' =  A_N x_N + B_{1N} w_N + B_{2N} u_N \\
z_N  =  C_{1N} x_N + D_{12N} u_N \\
y_N  =  C_{2N} x_N + D_{21N} w_N
\end{array} \right.
\end{equation}
where the operators of system (\ref{P}) have been replaced by
real-valued matrices computed on truncated hermitian basis.
We denote by $N_B$ the number of eigenfunctions of operator $\Delta^2$ we consider
and by $N_Z$ the number of Zernike modes used to describe $\phi_\tn{tur} $.
Then, $x_N\in\mathbb{R}^{2N_B + N_Z}$ is the state vector,
$w_N\in\mathbb{R}^{2N_B + 2N_Z}$ is the exogenous perturbation vector,
$u_N\in\mathbb{R}^{N_B}$ is the control vector,
$z_N\in\mathbb{R}^{N_B+N_Z}$ is the controlled output vector and
$y_N\in\mathbb{R}^{N_B+N_Z}$ is the measured output vector. The
matrices $A_N$, $B_{1N}$, $B_{2N}$, $C_{1N}$,
$D_{12N}$, $C_{2N}$ and $D_{21N}$ are of appropriate dimensions.

In order to compute these objects, we still consider the case of a
circular bimorph mirror which is free at all the boundary (this is also
the case of the mirror considered in Section \ref{numerical} below). The eigenvectors of operator
$$-\frac{Q_1}\rho \Delta^2 - \frac{Q_2}\rho I$$
are given by, for
all $ (k,j)\in \mathbb{N}^2$,
\begin{eqnarray*}
&L_{kj}(r,\theta)= a_{kj}\left(J_k \left(\dfrac{\lambda_{kj} r}{a }\right)
+ c_{kj} I_k \left(\dfrac{\lambda_{kj} r}{a }\right) \right)\cos(k\theta) \\
&
M_{kj}(r,\theta)= a_{kj}\left(J_k \left(\dfrac{\lambda_{kj} r}{a }\right)
+ c_{kj} I_k \left(\dfrac{\lambda_{kj} r}{a }\right) \right)\sin(k\theta)
\end{eqnarray*}
where $(r,\theta)$ are the polar coordinates of $x\in\Omega$, $J_k$ and $I_k$ are, respectively,
ordinary and modified Bessel function of first kind and order $k$, and
$ -\frac{Q_1}\rho \left(\frac {\lambda_{kj}}a\right)^4 - \frac{Q_2}\rho$ the corresponding eigenvalues.
The family $$\left\{  L_{kj}, M_{kj}, (k,j)\in \mathbb{N}^2  \right\}$$ is an Hilbertian basis of $H^2_{bc}(\Omega)$.
The dimensionless coefficients $\lambda_{kj}$  and $c_{kj}$ depend on the boundary conditions
while $a_{kj}$is computed using a normalization condition on the eigenvectors (see \cite{A-P-D} for further details).
In what follows, we consider the case of Poisson ratio $\nu=0.2$ corresponding to the material the mirror is made of.
Once a maximal azimuthal order is given (here $k_{\textnormal{max}}= 5$) the modes are classified
according to increasing $\lambda_{kj}$ and one has the values gathered in Table \ref{lambda}.

\begin{table}
\begin{center}
\begin{tabular}{|ccc|c|c|c|}
\hline
$i$&$j$&$k$&\(\lambda_{kj}\)&\(c_{kj}\)&\(a_{kj}\)\\
\hline
1&0&2&2.37805&0.18773&3.6157\\
2&1&0&2.96173&-0.092478&2.1984\\
3&0&3&3.60924&0.075982&4.4749\\
4&1&1&4.51025&-0.019949&3.8317\\
5&0&4&4.76934&0.034281&5.2453\\
6&0&5&5.89565&0.016333&5.9506\\
7&1&2&5.94302&-0.0056226&4.4178\\
8&0&2&6.18269&0.0032602&3.1394\\
9&1&3&7.30051&-0.0018233&4.9425\\
10&2&1&7.72338&0.0007269&4.9616\\
\hline
\end{tabular}
\caption{Coefficients of the eigenvectors $L_{kj}$ and $M_{kj}$.}\label{lambda}
\end{center}
\end{table}

The sequence of functions $L_{kj}$ and $M_{kj}$ need to be re-ordered. They are now denoted by $B_n$ and follow
the increasing values of $\lambda_{kj}$, alternating cosine and sine and eliminating the null eigenvectors $M_{0j}$.
Therefore,
$$ \forall x\in X,~ x = \sum_{n\in \mathbb{N}, \; i \geq 1} \alpha_{i} B_i(r,\theta)$$
where $(\alpha_{n})_{n\geq 1}$ is a sequence of real numbers satisfying $\sum_{n\in \mathbb{N}, \; n\geq 1} \alpha_{n}^2 <\infty$.\\

In reference \cite{B}, one can find that this basis $(B_n)_{n\in \mathbb{N}}$ with free boundary conditions is not orthogonal
in $L^2(\Omega)$. However, numerically, we can prove that this basis is nearly orthogonal, indeed lots of scalar products
in $L^2(\Omega)$ are null and the others are small $(10^{-6})$ in comparison with unity. So, for more numerical facilities,
we will use the scalar product in $L^2(\Omega)$ rather than in  $H^2_{bc}(\Omega)$.

Given $N_B$ and $N_Z \in \mathbb{N}$, we compute $A_N$, $B_{1N}$, $B_{2N}$, $C_{1N}$,
$C_{2N}$, $D_{12N}$ and $D_{21N}$ using the ``Bessel" truncated basis
$\{ B_{0}, B_{1}, \ldots, B_{N_B}\}$ and the Zernike one $\{ Z_{0}, Z_{1}, \ldots, Z_{N_Z}\}$.

We make analogous assumptions for the tuning parameters $b$, $c$ and $d$, i.e. $b=\mathrm{diag}_{i} (b_i)$, $c=\mathrm{diag}_{i} (c_i)$ and $d=\mathrm{diag}_{i} ( d_i)$ where $(b_i)_{i\in\mathbb{N}, \;i\geq 1}$, $(c_i)_{i\in\mathbb{N}, \;i\geq 1}$ and $(d_i)_{i\in \mathbb{N}, \; i\geq 1}$ are sequences of real numbers. We recall that these coefficients are weighting functions defining the respective weights of the disturbance signals and the choice of diagonal matrices corresponds to an assumption of decoupling between the different modes.

Futhermore $\phi_\tn{res}$ is expressed on Besssel functions, so we need to estimate
a projection matrix to define $\phi_\tn{tur}$ with Bessel spatial coordinates. We note $Q$ this projection
$N_B\times N_Z$-dimension matrix. Thus, the computed equation becomes:
\[
\phi_{\tn{res},i}=-\dfrac{4\pi}{\lambda}e_i+\sum_{j=1}^{N_Z-2}Q_{ij}B_{j+2}
\]
We denote by $\bf{0}$ each null matrix with the appropriate dimensions so that each following matrix makes sense. We get
\[
A_N=
\left [\begin{array}{ccc}
\bf 0 & {\bf 1}_{N_B} & \bf 0\\
- \omega_i^2 {\bf 1}_{N_B} & \bf 0 &\bf 0\\
\bf 0 & \bf 0 & F
\end{array}\right ]~~
B_{1N}=
\left [\begin{array}{cccc}
\bf 0 & \bf 0 & \bf 0 & \bf 0\\
 b & \bf 0 & \bf 0 & \bf 0\\
\bf 0 & \bf 0 & G & \bf 0
\end{array}\right ]\\
\]

\[
B_{2N}=\left [\begin{array}{cc}
\bf 0\\
\mathrm{block}_{ij}\left (\dfrac{\tilde{d}_{31}}{\rho}\left<\Delta B_i, B_j \right> \right ) \\
\bf 0
\end{array}\right ]~~
C_{1N}=\left [\begin{array}{ccc}
-\dfrac{4\pi }{\lambda}{\bf 1}_{N_B} & \bf 0& Q \\
 \bf 0 & \bf 0 & \bf 0
 \end{array}\right ]~~
 D_{12N}=\left [\begin{array}{cc}
\bf 0\\
{\bf 1}_{N_B}
\end{array}\right ]
\]

\[
C_{2N}=\left [\begin{array}{cccc}
\mathrm{block}_{ij}\left (\tilde{e}_{31} \left<\Delta B_i, B_j \right> \right )& \bf 0 &\bf 0\\
-\dfrac{4\pi}{\lambda}{\bf 1}_{N_B} & \bf 0 &Q
 \end{array} \right ]~~
D_{21N}=\left [\begin{array}{cccccccc}
\bf 0 &\bf 0&\bf 0&  d\\
\bf 0 & c &\bf 0 &\bf 0\\
\end{array}\right ]
\]
where $\omega_i^2 = \dfrac{Q_1}{\rho} \left(
\dfrac{\lambda_{i}}{a} \right)^4+ \dfrac{Q_2}{\rho}$
and $\left< \cdot ,  \cdot \right>$ is the usual scalar product in $L^2(\Omega)$.

%%%%%%%%%%%%%%%%%%%%%%%%%%%%%%%%%%%%%%%%%%%%%%%%%
\subsection{Numerical results}\label{numerical}
%%%%%%%%%%%%%%%%%%%%%%%%%%%%%%%%%%%%%%%%%%%%%%%%%

In this subsection, numerical simulations are proposed.
To get more realistic results, the experimental device of the project SESAME of the Observatoire de
Paris is considered. This experimentation uses a bimorph mirror with a
distribution of 31 piezoelectric actuators. The piezoelectric inclusions are PZT patches. We use the physical constants of Table \ref{PC}

\begin{table}
\begin{center}
\begin{tabular}{|c|c|}
\hline
Wind speed &$V = 9$ ms$^{-1}$\\
\hline
Diameter of the pupil&$D = 10^{-2}$ m\\
\hline
radius of the mirror&$a = 25\times10^{-3}$ m\\
\hline
mirror's stiffness coefficients&$Q_1 = 84$ Nm,~   $Q_2 = 11.25 \times10^8$ Nm$^{-3}$\\
\hline
mirror's surfacic density & $\rho = 16.3$ kg.m$^{-2}$\\
\hline
piezoelectric coefficients & $\tilde d_{31} = -0.0044$ NV$^{-1}$,~  $ \tilde{e}_{31}=-5.60\times 10^{3}$ Vm \\
\hline
wave length & $\lambda=550$ nm\\
\hline
\end{tabular}
\vspace{0.2cm}
\caption{Physical parameters for the numerical simulations}\label{PC}
\end{center}
\end{table}

\begin{figure}[h]
\centerline{\includegraphics[width=9cm]{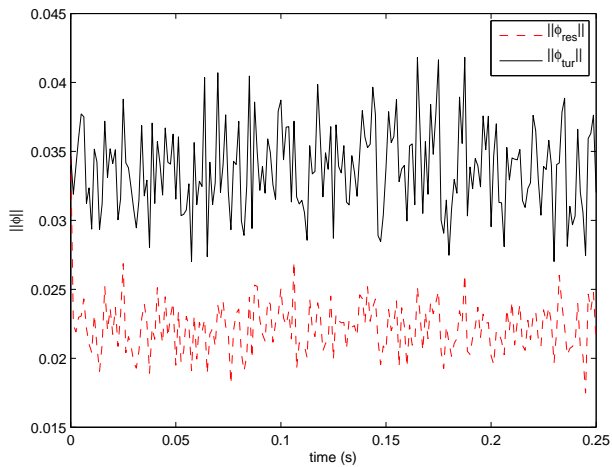}}
\caption{Time-evolution of $\|\phi_\tn{tur}\|_{L^2(\Omega)}$ (solid line) and $\|\phi_\tn{res}\|_{L^2(\Omega)}$ (dashed line)}\label{Turbulence}
\end{figure}

We simulate only the $12$ modes which follow the tip/tilt.
The performance of the control system is evaluated by considering the spatial norm $\|.\|_{L^2}$ of $\phi_\tn{res}$
compared to $\|\phi_\tn{tur}\|_{L^2(\Omega)}$:
\[
\|\phi_\tn{tur}\|_{L^2(\Omega)}=\sum_{i=4}^{N_z}\phi_i(t)^2.
\]
For identical random initial conditions and taking the respective weights of the disturbance signals such that $b_i=0.001$, $c_i=0.002$ and
$d_i=0.003$ for all $i$, we obtain the results represented in Figure \ref{Turbulence}.

Using Monte Carlo simulations, the ratio between temporal average of  $\|\phi_\tn{tur}\|_{L^2(\Omega)}$ and $\|\phi_\tn{res}\|_{L^2(\Omega)}$ is near
to 1.91 which represents a phase distortion attenuation of the reflected wavefront of $48 \%$. In addition one should recall that this result does
not take into account the tip/tilt correction. Even if these results are of the same order of magnitude as those presented in \cite{PAS93}, which
cannot be considered as completely satisfactory when considering usual results on real experiments, they clearly demonstrate the feasibility of the
proposed approach. The possible degradation of such a performance induced by the delay in the loop and the discretization of the control law for
implementation purpose could darken the picture. It must be recalled that this apparent loss of performance is mainly due to the tuning of the
trade-off between robustness and performance that is inherently encountered in closed-loop feedback design. Numerous improvements have still to be
considered as presented in the next conclusion.

\section{Conclusion}
In this paper, a new framework to deal with the problem of adaptive optics is proposed. It is mainly based on an infinite-dimensional model of the
deformable mirror associated with the definition of a standard model on which robust control techniques may be applied. The preliminary numerical
experiments show a performance level comparable with the results of reference \cite{PAS93}. The main advantage of the approach suggested in this
paper is that no interaction matrix is required to control the system. We do not pretend to outperform already existing AO systems but rather to pave
the way for future major improvements in terms of robustness and efficiency of the proposed control strategies. The authors are planing to take
into account a model for the Shack-Hartmann wavefront sensor including a time delay associated with processing measurements. This will be covered in
a next study.

\section*{Acknowledgements} The authors are grateful to Pascal Jagourel, Observatoire de Paris, for useful discussions on bimorph mirrors.

\end{document}